\documentclass{amsart}
\usepackage{amssymb, latexsym}

\begin{document}
\title{On Closed Subsets of Free Groups}
\author{Rita Gitik}
\address{ Department of Mathematics \\ University of Michigan \\ Ann Arbor, MI, 48109}
\email{ritagtk@umich.edu}
\author{Eliyahu Rips}
\address{ Institute of Mathematics \\ Hebrew University, Jerusalem, 91904, Israel}
\date{\today}

\begin{abstract}
We give two examples  of a finitely generated subgroup of a free group 
and a subset, closed in the profinite topology of a free group, such that their product
is not closed in the profinite topology of a free group.
\end{abstract}


\maketitle


\section{Introduction}

A theorem of M. Hall, proved in \cite{Hall}, states  
that any finitely generated subgroup of a free group is closed 
in the profinite topology. 
This result has been  generalized by many researchers.
  
The authors proved in \cite{G-R1} and \cite{G-R2} 
that the product of two finitely generated subgroups of a free group
is closed in the profinite topology of the free group. 
The first published proof of that theorem is due to G.A. Niblo, cf. \cite{Ni}.
 
Denote the profinite topology on a free group $F$ by $PT(F)$.
 
A more general result saying that for any finitely generated subgroups
$H_1, \cdots ,H_n$ of a free group $F$  the set 
$H_1 \cdots H_n$ is closed in  $PT(F)$ was obtained by L. Ribes and P.A. Zalesskii in
\cite{R-Z}, by K. Henckell, S.T. Margolis, J.E. Pin, and J. Rhodes in \cite{H-M-P-R},
and by B. Steinberg in \cite{Ste}.
 
T. Coulbois in \cite{Co} proved that property $RZ_n$ is closed under free products, 
where a group $G$ is said to have property $RZ_n$ if for any $n$ finitely generated subgroups 
$H_1, \cdots  H_n$ of $G$, the set $H_1 \cdots H_n$ is closed in $PT(G)$. 

The aforementioned results lead to the following question:
is it true that for any finitely generated subgroup $H$ of a free
group $F$ and for any subset $S$ of $F$ which is closed in $PT(F)$, the
product $SH$ is closed in  $PT(F)$.

In this paper we provide a negative  answer to this question by constructing
two counterexamples.

\section{The First Example}

Let $F$ be a finitely generated free group. The profinite topology on $F$ is defined by
proclaiming all subgroups of finite index of $F$ and their cosets to be basic open sets.
An open set in $PT(F)$ is a (possibly infinite) union of 
cosets of various subgroups of finite index and
the closed sets in $PT(F)$ are the complements of such unions in $F$.

The following example describes a set $S$, closed in $PT(F)$, such that its product with 
a free factor of $F$ is not closed in $PT(F)$.

\textbf{Example 1.}

Let $F = <a,b>$ be a free group of rank two.
Consider an infinite sequence $A = \{ a, a^{2!}, a^{3!}, \cdots, a^{k!}, \cdots \} \subset F$.
Note that $A$ converges to $1_F$. Indeed, let $N$ be a normal subgroup of finite index $m$ in $F$.
If $k \ge m$  then $a^{k!}$ is contained in $N$. 
Hence any open neighborhood of $1_F$ in $F$ contains all, but finitely many elements of $A$,
therefore $A$ converges to $1_F$. Note that $1_F \notin A$, so $A$ is not closed in $PT(F)$. 

Let $m_k, k \ge 1$ be integers such that 
$m_k \rightarrow m_0 \in \hat{Z} \setminus Z$ in $PT(\hat{Z})$,
where $\hat{Z}$ is the completion of $Z$ in $PT(Z)$.
Then $a^k b^{m_k} \rightarrow a^0 b^{m_0} \in \hat{F} \setminus F$,
where $\hat{F}$ is the completion of $F$ in $PT(F)$.
Hence the sequence $a^k b^{m_k}$ has no other limit points.
In particular, it has no limit points in $F$. 
Therefore for every $w \in F$ with 
$w \neq a^k b^{m_k}$ for all $k \ge 1$, there exists an open neighborhood $U$ of $w$
such that $a^k b^{m_k} \notin U$, for all $k \ge 1$.
It follows that the set 
$S= \{ a b ^{m_1}, a^{2!} b ^{m_2}, \cdots , a^{k!} b ^{m_k}, \cdots \}$ is closed in $PT(F)$.

Note that $1_F \notin S<b>$, however $ A \subseteq S<b>$, so
$1_F \in \bar A \subseteq \overline{S<b>}$.
We conclude that $S<b>$ is not closed in $PT(F)$.

\section{The Second Example}

The example in the previous section raises the following question:
is it possible to impose some restrictions on a set $S$, closed in $PT(F)$,
such that the product of $S$ with a free factor of $F$ would be closed in $PT(F)$.

The following example demonstrates that such restrictions on $S$ should be severe.

Let $F$ be a free group on free generators $K \cup L$, with $|K|=k$, $|L|=l$, and
$F = <K>*<L>$. We describe a discrete set $S$, closed in $PT(F)$, such that $S<K>$ 
is not closed in $PT(F)$ and the last syllable of all elements of $S$ is in $<L>$.

\textbf{Example 2.}

Construct by induction a sequence of normal subgroups of finite index

$G_1 > G_2 > \cdots > G_m > \cdots $, 
elements $r_m \in <K>$ and $s_m \in F$, and an increasing function $f(m)$
satisfying the following conditions:

\begin{enumerate}
\item
$G_m \cdot r_m$ is at distance greater than $f(m)$ from $G_m \cdot 1$ in $F/G_m$.
\item
$G_ms_m = G_m r_m$.
\item
The last syllable of all $s_m$ is in $<L>$.
\item
For all $k>m, G_mr_k \neq G_mr_m$.
\end{enumerate}

Let $f(1)$ be an arbitrary integer. Choose $G_1$ 
such that the index of 

$H_1=G_1 \cap <K>$ in $<K>$ exceeds $2k(2k-1)^{f(1)-1}$, which 
is the upper bound on the number of elements in a ball of radius $f(1)$ in $<K>/H_1$  around
$1$. Then we can choose $r_1 \in <K>$ such that the distance between $G_1 \cdot r_1$ and $G_1 \cdot1$
in $F/G_1$ is bigger than $f(1)$. 
Choose $s_1 \in F$ such that $G_1r_1=G_1s_1$ and the last syllable of $s_1$ is in $<L>$.

Assume that for some $n > 1$ we have constructed the normal subgroups

$G_1 > G_2> \cdots > G_n$ of finite index in $F$, elements $r_1, \cdots , r_n $ of $<K>$
and $s_1, \cdots, s_n$ of $F$, and $f(1) < f(2) < \cdots < f(n)$ satisfying conditions
$1$, $2$, $3$, and $4$ for any $m \le n$.

Let $H_m = G_m \cap <K>, m=1, \cdots , n$.
Let $\phi_{i,j} : F/G_i \rightarrow F/G_j$ for $i > j$ be the natural homomorphisms, and let
$\psi_{i,j}:<K>/H_i \rightarrow <K>/H_j$ be the corresponding natural homomorphisms.

We want to define $G_{n+1}$, $r_{n+1} \in <K>$, $s_{n+1} \in F$, and $f(n+1)$.

In order to satisfy condition $4$, we need $G_{n+1} r_{n+1}$  to be distinct from the cosets

 $\phi^{-1}_{n+1,1}(G_1r_1), \phi^{-1}_{n+1,2} (G_2r_2), \cdots , \phi^{-1}_{n+1,n}(G_n r_n)$.

In order to choose $G_{n+1}r_{n+1}$ such that $r_{n+1} \in <K>$, we need to
consider the subgroups $H_i = G_i \cap <K>$ and the preimages

$\psi^{-1}_{n+1,1}(H_1r_1)$, $\psi^{-1}_{n+1,2}(H_2r_2), \cdots , \psi^{-1}_{n+1,n}(H_nr_n)$.

Note that the preimage of $H_mr_m$ in $<K>/H_{n+1}$ contains $[H_m : H_{n+1}]$ elements. 
Therefore, $H_{n+1}r_{n+1}$ should be different from
$[H_1 : H_{n+1} ] +[H_2 : H_{n+1} ]+ \cdots +[H_n : H_{n+1} ]$ out of the total $[<K> : H_{n+1} ]$
elements. In addition, in order to satisfy condition $1$, the coset $H_{n+1}r_{n+1}$ should lie 
outside the ball of radius $f(n+1)$ around $1$ in $<K>/H_{n+1}$. 
Note that the number of elements in this ball does not exceed $2k(2k-1)^{f(n+1)-1}$.

Also note that $[H_1 : H_{n+1} ] +[H_2 : H_{n+1} ]+ \cdots +[H_n : H_{n+1} ] =$

$=[<K> : H_{n+1}]({1 \over [<K> : H_1]} + {1 \over [<K> : H_2]} + \cdots + 
{1 \over [<K> : H_n]})$.

Assuming that the sequence of indices $[<K> : H_1]$, $[<K> : H_2], \cdots, $
increases rapidly enough, we may assume that for all $n \ge 1$ the quantity

${1 \over [<K> : H_1]} + {1 \over [<K> : H_2]} + \cdots + 
{1 \over [<K> : H_n]}$ is smaller than $1 \over 2$.

Now we can choose $G_{n+1}$ such that the index  $[F : G_{n+1}]$ is big enough and for
 
$H_{n+1} = <K> \cap G_{n+1}$ the index $[<K> : H_{n+1}]$ is big enough.
 
We choose an element $r_{n+1} \in <K>$ such that $G_{n+1} r_{n+1}$ is outside the ball
of radius $f(n+1)$ in $F/G_{n+1}$ around $G_{n+1} \cdot 1$ and such that
$G_{n+1}r_{n+1}$ is different from all the elements of the preimages
$\phi^{-1}_{n+1,1}(G_1r_1), \phi^{-1}_{n+1,2} (G_2r_2), \cdots , \phi^{-1}_{n+1,n}(G_n r_n)$.

 Choose $s_{n+1} \in F$ such that $G_{n+1}s_{n+1} = G_{n+1}r_{n+1}$ and the last syllable of $s_{n+1}$
 is in $<L>$.
 
Hence, by induction, we have satisfied conditions $1$, $2$, $3$, and $4$.
 
Let $S = \{s_1, s_2, \cdots \}$.
We claim that $S$ is discrete in $PT(F)$. Indeed, for each $n$ the coset $G_n r_n =G_n s_n$
does not contain any $G_nr_k =G_ns_k$ for $k > n$, so we have found an open neighborhood
 of $s_n$ containing at most $n$ elements of $S$. As the profinite topology on a free
 group is Hausdorff, it follows that $S$ is discrete. 
 
Note that $S$ does not have limit points in $F$. Indeed, consider $x \in F$.
For any $n \ge 1$ the coset $G_n \cdot x$ is an open neighborhood of $x$ in $PT(F)$
and the distance between $G_n \cdot 1$ and $G_n \cdot x$ is bounded by the length of $x$.

By definition of $S$, the distance between $G_n \cdot 1$ and $G_n \cdot s_n$ is greater
than $f(n)$ for almost all $n \ge 1$, hence  the intersection $S \cap G_n \cdot x$ is finite.

It follows that $x$ is not a limit of $S$, therefore, $S$ is closed in $PT(F)$.

Note that $1 \in G_ns_nr^{-1}_n \in G_nS<K>$ for all $n \ge 1$, so
$1 \in \overline{S<K>}$, but $1 \notin S<K>$ because $S \cap <K> = \emptyset$.
Therefore, $S<K>$ is not closed in $PT(F)$.
  
\section{Acknowledgment}

The first author would like to thank the Institute of Mathematics
of the Hebrew University for generous support.


\begin{thebibliography}{10}
 
\bibitem{Co}
T. Colulbois,
\emph{Free Product, Profinite Topology, and Finitely Generated Subgroups},
Int. J. of Algebra and Computation,
\textbf{11} (2001), 171-184.
 
\bibitem{G-R1} 
R. Gitik and E. Rips,
\emph{On Separability Properties of Groups},
Int. J. of Algebra and Computation,
\textbf{5} (1995), 703-717.

\bibitem{G-R2}
R. Gitik and E. Rips,
\emph{On Double Cosets in Free Groups},
J. of Mathematics, IOSR,
\textbf{12}, Issue 6 Version IV (November-December 2016), 14-15.
  
\bibitem{Hall} 
M. Hall, Jr.,
\emph{Coset Representations in Free Groups},
Trans. AMS,
\textbf{67} (1949), 431-451.

\bibitem{H-M-P-R} 
K. Henckell, S.T. Margolis, J.E. Pin, and J. Rhodes, 
\emph{Ash's Type II Theorem, Profinite Topology, and Malcev Products I},
Int J. of Algebra and Computation,
\textbf{1} (1991), 411-436.

\bibitem{Ni}
G. A. Niblo, 
\emph{Separability Properties of Free Groups and Surface Groups}, 
J. of Pure and Applied Algebra,
\textbf{78} (1992), 77-84.

\bibitem{R-Z} 
L. Ribes and P.A. Zalesskii, 
\emph{On the Profinite Topology on a Free Group I}, 
Bull. LMS,
\textbf{25} (1993), 37-43.

\bibitem{Ste}
B. Steinberg,
\emph{Inverse Automata and Profinite Topologies on a Free Group},
J. of Pure and Applied Algebra,
\textbf{167} (2002), 341-359.


\end{thebibliography}
\end{document}